\documentclass{article}
\usepackage{ijcai20}

\usepackage{times}

\usepackage{soul}
\usepackage{url}
\usepackage[hidelinks]{hyperref}
\usepackage[utf8]{inputenc}
\usepackage[small]{caption}
\usepackage{graphicx}
\usepackage{amsmath}
\usepackage{booktabs}
\usepackage{pdfpages}
\usepackage{amsmath,amssymb,amsfonts,amsthm}       
\usepackage{mathtools}
\usepackage{xspace}
\usepackage{microtype}
\usepackage{todonotes}
\usepackage{colonequals}
\usetikzlibrary{patterns,arrows,arrows.meta,calc,shapes,shadows,decorations.pathmorphing,decorations.pathreplacing,automata,shapes.multipart,positioning,shapes.geometric,fit,circuits,trees,shapes.gates.logic.US,fit}
\usetikzlibrary{backgrounds,scopes}
 \usepackage{pdfpages}
\usepackage{multirow}
\usepackage{multicol}

\usepackage{paralist}

\usepackage{subfigure}
\usepackage{url}
\usepackage{tikz} 
\usepackage{mdframed}

\usetikzlibrary{arrows,decorations.pathmorphing,positioning,fit,trees,shapes,shadows,automata,calc} 

\theoremstyle{definition}
\newtheorem{definition}{Definition}

\theoremstyle{remark}
\newtheorem{remark}{Remark}
\theoremstyle{problem}

\newtheorem*{problem*}{Problem}

\title{Robust Policy Synthesis for Uncertain POMDPs via Convex Optimization}
 \author{
Marnix Suilen,\textsuperscript{\rm 1}
Nils Jansen,\textsuperscript{\rm 1}
Murat Cubuktepe,\textsuperscript{\rm 2} 
Ufuk Topcu,\textsuperscript{\rm 2} \\
\textsuperscript{\rm 1} Department of Software Science, Radboud University, The Netherlands\\
\textsuperscript{\rm 2} Department of Aerospace Engineering and Engineering Mechanics, University of Texas at Austin, USA
}

\newcommand{\ie}{i.e.\@\xspace}

\newcommand{\tool}[1]{\textrm{#1}\xspace}

\newcommand{\p}{\ensuremath{\mathbb{P}}}

\newcommand{\reachProp}[2]{\ensuremath{\p_{\geq #1}(\finally #2)}}
\newcommand{\reachProplT}{\ensuremath{\reachProp{\lambda}{T}}}

\newcommand{\reachPropSymbol}{\varphi_r}

\newcommand{\ereachPropSymbol}{\varphi_c}

\newcommand{\expRewProp}[2]{\ensuremath{\EV_{\leq #1}(\finally #2)}}

\newcommand{\finally}{\lozenge}

\newcommand{\R}{\mathbb{R}}

\newcommand{\Ireal}{[0,\, 1]\subseteq\mathbb{R}}  
\newcommand{\Ex}{\ensuremath{\mathbb{E}}\xspace}        
\newcommand{\EV}{\ensuremath{\mathbb{E}}}

\newcommand{\Distr}{\mathit{Distr}}

\newcommand{\distDom}{X}

\newcommand{\distFunc}{\mu}
\newcommand{\distDomElem}{x}

\newcommand{\sinit}{s_{\mathit{I}}} 
\newcommand{\umdp}{M}
\newcommand{\umdpInit}[1][]{\ensuremath{\umdp{#1}=(S{#1},\sinit{#1},\Act,\probumdp{#1},\intervals)}}

\newcommand{\probumdp}{\mathcal{P}}
\newcommand{\intervals}{\mathbb{I}}

\newcommand{\sched}{\ensuremath{\sigma}}

\newcommand{\Act}{\ensuremath{\mathit{Act}}}

\newcommand{\act}{\ensuremath{\alpha}}

\DeclareMathAlphabet{\mathpzc}{OT1}{pzc}{m}{it}
\def\presuper#1#2%
  {\mathop{}%
   \mathopen{\vphantom{#2}}^{#1}%
   \kern-\scriptspace%
   #2}

\renewcommand{\paragraph}[1]{\par\smallskip\noindent\textbf{#1}}

\newcommand{\mdp}{M}

\newcommand{\ObsSym}{{Z}}
\newcommand{\ObsFun}{{O}}

\newcommand{\PomdpInit}[1][]{\pomdp{#1}=(\mdp{#1},\ObsSym{#1},\ObsFun{#1})}
\newcommand{\pomdp}{\mathcal{M}}

\newcommand{\states}{\ensuremath{S}}

\newcommand{\osched}{\ensuremath{\mathit{\sigma}}}
\newcommand{\oSched}{\ensuremath{\Sigma}}

\begin{document}

\maketitle

\begin{abstract}
We study the problem of policy synthesis for uncertain partially observable Markov decision processes (uPOMDPs).
The transition probability function of uPOMDPs is only known to belong to a so-called uncertainty set, for instance in the form of probability intervals.
Such a model arises when, for example, an agent operates under information limitation due to imperfect knowledge about the accuracy of its sensors.
The goal is to compute a policy for the agent that is robust against all possible probability distributions within the uncertainty set.
In particular, we are interested in a policy that robustly ensures the satisfaction of temporal logic and expected reward specifications.
We state the underlying optimization problem as a semi-infinite quadratically-constrained quadratic program (QCQP), which has finitely many variables and infinitely many constraints.
Since QCQPs are non-convex in general and practically infeasible to solve, we resort to the so-called convex-concave procedure to convexify the QCQP.
Even though convex, the resulting optimization problem still has infinitely many constraints and is NP-hard.
For uncertainty sets that form convex polytopes, we provide a transformation of the problem to a convex QCQP with finitely many constraints. 
We demonstrate the feasibility of our approach by means of several case studies that highlight typical bottlenecks for our problem. 
In particular, we show that we are able to solve benchmarks with hundreds of thousands of states, hundreds of different observations, and we investigate the effect of different levels of uncertainty in the models.
\end{abstract}

\section{Introduction}
\label{sec:introduction}

Partially observable Markov decision processes (POMDPs) model sequential decision-making problems under stochastic uncertainties and partial information~\cite{kaelbling1998planning}.
In particular, an agent that operates in an environment modeled by a POMDP receives observations according to which it tries to infer the likelihood, called the belief state, of the system being in a certain state.
Based on this partial information about the environment, the agent chooses action whose outcome is stochastically determined.


The assumption that the transition and observation probabilities in POMDPs are explicitly given does often not hold.
Unforeseeable events such as (unpredictable) structural damage to a system~\cite{796wwe3643} or an imprecise sensor model~\cite{bagnell2001solving}, may necessitate to account for additional uncertainties in the value of the probabilities.
So-called uncertain POMDPs (uPOMDPs) address this need by incorporating sets of uncertainties in the probabilities.
The sets may be described as, for example, intervals~\cite{givan2000bounded} or more generally by likelihood functions~\cite{DBLP:journals/ior/NilimG05}.
For example, take a robust aircraft collision avoidance system that issues advisories to pilots~\cite{kochenderfer2015decision}.
Modeled as a uPOMDP, the actions relate to such advice and concern the flying altitude and the speed.
Uncertainty enters the model in the form of unreliable data of the reaction time of a pilot.
Moreover, there is an uncertain probability of receiving false observations regarding the speed and altitude of other aircrafts.

We study the synthesis of policies in uPOMDPs.
Specifically, we seek to compute a policy that satisfies temporal logic~\cite{Pnueli77} or expected reward specifications against all possible probability distributions from the uncertainty set.
For the aforementioned collision avoidance system, such a policy would minimize the acceleration due to fuel efficiency and ensure that the probability of not colliding with another aircraft is above a certain threshold.

The \emph{robust synthesis problem} for uncertain MDPs, that is, with full observability, has been extensively studied.
The existing approaches rely, for instance, on dynamic programming~\cite{DBLP:conf/cdc/WolffTM12}, convex optimization~\cite{DBLP:conf/cav/PuggelliLSS13}, or value iteration~\cite{DBLP:conf/qest/HahnHHLT17}.
The existing approaches for uPOMDPs, rely, for example, on sampling~\cite{burns2007sampling} or robust value iteration~\cite{DBLP:conf/icml/Osogami15} on the belief space of the uPOMDP, but do not take into account temporal logic constraints.
In general, the robust synthesis problem for uPOMDPs is hard. 
Computing an optimal policy even for POMDPs, that is, with no uncertainty in the probabilities, is undecidable in general~\cite{MadaniHC99}, PSPACE-complete~\cite{meuleau1999learning} for finite-horizon properties, and NP-hard~\cite{VlassisLB12} if the polices do not take the execution history into account, that is, they are memoryless.

We develop a novel solution for efficiently computing policies for uPOMDPs using robust convex optimization.
We restrict the problem to memoryless policies, while finite memory may be added effectively~\cite{junges2018finite}.
Moreover, we focus---for brevity in this paper---on uncertainty sets that are given by intervals, \ie, upper and lower bounds on probabilities. 
The approach, though, is applicable to all uncertainty sets that are represented as convex polytopes.

First, we encode the problem as a semi-infinite quadratically-constrained quadratic program (QCQP), which includes finitely many variables but infinitely (in fact uncountably) many constraints that capture the uncertainty set~\cite{wiesemann2013robust}.
The structure of the encoding is similar to the one for POMDPs without uncertainty~\cite{amato2006solving}.
This optimization problem is non-convex in general and thereby infeasible to solve in practice~\cite{DBLP:conf/nips/ChenLSS17}.
We use the so-called convex-concave procedure to convexify the problem~\cite{lipp2016variations}.
The resulting convex QCQP provides a sound over-approximation of the original problem, yet, it still has infinitely many constraints and renders the application of the convex-concave procedure impractical in that form. 

Towards computational tractability for solving the semi-infinite convex QCQP, we restrict the uncertainty set to convex polytopes and gain two key advantages.
First, a convex polytope represents the \emph{valid probability distributions} exactly and avoids an unnecessarily coarse approximation of the uncertainty.
Second, it suffices to enumerate over the finite set of vertices of these polytopes to retain optimal solutions~\cite[Section 5.2]{Lofberg2012automaticrobust}.
We exploit this property and transform the semi-infinite program to a finite convex QCQP which, integrated with the convex-concave procedure, provides an efficient solution to the robust synthesis problem.

Three complicating factors require special attention in the proposed solution.
First, the iterative convex-concave procedure (CCP) may take exponentially many iterations in the number of its input variables~\cite{park2017general}. 
The reason is that the standard stopping criterion of the CCP is conservative, and we observe in the numerical examples that it largely affects the runtime.
We provide a dedicated version of the CCP that mitigates this problem by integrating a robust verification method~\cite{benedikt2013ltl}, similar as in~\cite{cubuktepe-et-al-qcqp-techreport} for so-called parametric MDPs~\cite{DBLP:journals/corr/abs-1903-07993}.
In particular, we compute the exact probability and expected cost values in intermediate candidate solutions delivered by the CCP.
We check whether these solutions already satisfy the specifications, otherwise the exact values are used as input for the next iteration of the CCP.

Second, enumerating the vertices of the convex polytope causes an exponential number of constraints in the number of uncertain transitions.
This number, however, depends here on the number of successor states of each state-action pair in the uPOMDP, and typical benchmarks, as available at \texttt{\url{http://pomdp.org}} or used in~\cite{NPZ17}, are usually sparse, reducing the effect of this theoretical blowup.

The third complicating factor is the general hardness of problems with partial observability and particularly due to the number of observations. 
On the other hand, the size of the resulting convex optimization problems in the proposed solution is polynomial in the number of observations as well as states.
Note that the range of the uncertainty sets, or more specifically the size of the intervals, does not affect the number of constraints.
With our prototype implementation, we solve problems with hundreds of thousands of states and thousands of observations for several well-known case studies.

\paragraph{Related work.}
To the best of our knowledge, the proposed approach is the first that accounts for temporal logic specifications in the computation of policies for uPOMDPs. 
Beyond that,~\cite{burns2007sampling} relies on sampling and~\cite{DBLP:conf/icml/Osogami15} uses robust value iteration on the belief space of the uPOMDP. 
\cite{itoh2007} assumes a belief over the uncertainty, that is, distributions over the exact probability values.
Robustness in~\cite{DBLP:conf/cdc/ChamieM18} is defined over fixed belief regions. 

In addition, one could adapt the approaches that employ mixed-integer linear programming for POMDPs~\cite{DBLP:conf/aips/ArasDC07,DBLP:conf/aips/KumarMZ16} by defining upper and lower bounds on the probabilities.
 In that case, a solution to the optimization problem would induce a rather coarse over-approximation because it will necessarily take into account sub- and super-distributions. 
Contrarily, we compute tight approximations of optimal robust policies and demonstrate this in our numerical examples.


\section{Preliminaries}
\label{sec:preliminaries}

A \emph{probability distribution} over a finite or countably infinite set $\distDom$ is a function $\distFunc\colon\distDom\rightarrow\Ireal$ with $\sum_{\distDomElem\in\distDom}\distFunc(\distDomElem)=1$. 
The set of all distributions on $\distDom$ is denoted by $\Distr(\distDom)$.
%
	A convex polytope is a convex $n$-dimensional shape defined by $n$ linear inequations of the form $A\vec{x} \leq \vec{c}$, where $A \in \mathbb{R}^{n \times m}$ and $\vec{c} \in \mathbb{R}^{n}$.



\begin{definition}[uMDP]
\label{def:umdp}
  An \emph{uncertain Markov decision process (uMDP)} is a tuple $\umdpInit$ where $\states$ is a set of states, $\sinit \in \states$ is the initial state, $\Act$ is the set of actions, $\intervals = \{[a,b] \mid a,b \in (0,1] \text{ and } a \leq b\}$ is a set of probability intervals, such that $\probumdp \colon \states \times \Act \times \states \to \intervals$ forms the \emph{uncertain transition function}. 
  A reward function $r \colon \states \times \Act \to \R_{\geq 0}$ assigns rewards to state action pairs.
\end{definition}

%
%
%
%
For a uMDP $\mdp$ and a transition probability function $P \colon \states \times \Act \to \Distr(\states)$, we write $P \in \probumdp$ if for all $s,s' \in \states$ and $\alpha \in \Act$ we have $P(s,\act,s') \in \probumdp(s,\act,s')$. 
Intuitively, $P$ yields only values within the corresponding intervals of $\probumdp$ for each state-action pair $s,\alpha$ and corresponding successor state $s'$.
We restrict the function $P$ to only select values from the intervals that form valid probability distributions and discuss later how this is achieved algorithmically.
A uMDP is \emph{instantiated} by $P \in \probumdp$, yielding a Markov decision process (MDP) $\mdp[P]$.

\begin{remark}
	For the correctness of our method, we require lower bounds of intervals to be strictly larger than $0$. 
An instantiation can not ``eliminate'' transitions from the uMDP by assigning value $0$.
The problem statement would be different and theoretically harder to solve if we also include values of $0$ for the transitions, see~\cite{winkler2019complexity}.
Furthermore, we allow the upper and lower bound of an interval to be the same, resulting in \emph{nominal} transition probabilities.
\end{remark}
%

%


\begin{definition}[uPOMDP]
   An \emph{uncertain partially observable MDP (uPOMDP)} is a tuple $\PomdpInit$, with $\umdpInit$ the \emph{underlying uMDP of $\pomdp$}, $\ObsSym$ a finite set of \emph{observations} and $\ObsFun\colon\states\rightarrow\ObsSym$ the \emph{observation function}.
\end{definition}
%
We assume that all states have the same actions. 
%
More general observation functions use a distribution
over $\ObsSym$, while there is a polynomial transformation of the general case to our definition~\cite{ChatterjeeCGK16}.
We define instantiations of uPOMDPs via the underlying uMDP.


%
  An \emph{observation-based policy} $\osched\colon Z\rightarrow\Distr(\Act)$ for a uPOMDP maps observations to distributions over actions. 
  Note that such a policy is referred to as memoryless and randomized. 
  More general (and powerful) types of policies take an (in)finite sequence of observations and actions into account.
  $\oSched^\pomdp$ is the set of observation-based strategies for $\pomdp$.
  Applying $\osched\in\oSched^\pomdp$ to $\pomdp$ resolves all choices and partial observability and an induced (uncertain) Markov chain $\pomdp^\osched$ results.
  
For a POMDP $\pomdp_e$ (without uncertainties) and a set of target states $T\subseteq \states$, the  \emph{reachability specification} $\reachPropSymbol=\reachProplT$ states that the probability of reaching $T$ shall be at least $\lambda$.
Similarly, the \emph{expected cost specification} $\ereachPropSymbol = \expRewProp{\kappa}{G}$ states that the expected cost of reaching the goal set $G \subseteq \states$ shall be less than or equal to $\kappa$.
A policy $\osched\in \oSched^{\pomdp_e}$ satisfies $\reachPropSymbol$ (or $\ereachPropSymbol$) if it is satisfied on the Markov chain $\pomdp^\osched$, denoted by $\osched\models\reachPropSymbol$ ($\osched\models\ereachPropSymbol$).
A policy for uPOMDPs takes all possible instantiations from the uncertainty sets into account. 
\begin{definition}[Robust Policy]
For a uPOMDP $\pomdp$, the underlying uMDP $\umdpInit$, and a specification $\varphi$, 
an observation-based policy $\osched\in\oSched^\pomdp$ \emph{robustly satisfies} $\varphi$ for $\pomdp$ ($\osched\models\varphi$) if for all $P\in\probumdp$ it holds that $\pomdp[P]^\osched$ satisfies $\varphi$.
\end{definition}
Intuitively, the policy needs to satisfy the specification for all instantiations from $\pomdp[P]$.
If we have several (expected cost or reachability) specifications $\varphi_1,\ldots,\varphi_m$, we write $\osched\models\varphi_1\land\ldots\land\varphi_n$ where $\osched$ robustly satisfies all specifications.
Note that general temporal logic constraints can be reduced to reachability specifications~\cite{BK08,DBLP:journals/corr/abs-2001-03809}, therefore we omit a detailed introduction to the underlying logic.

\section{Formal Problem and Outline}
We first state the central problem of this paper. 
\begin{mdframed}[backgroundcolor=gray!30]
\begin{problem*}[Robust Synthesis for uPOMDPs]\label{prob:upomdp}
Given an uPOMDP $\PomdpInit$ and a specification $\varphi$, which is either a reachability specification $\reachPropSymbol=\reachProplT$ or an expected cost specification $\ereachPropSymbol = \expRewProp{\kappa}{G}$, compute a randomized memoryless policy $\osched\in\oSched^\pomdp$ for $\pomdp$ such that $\osched$ robustly satisfies the specification, that is, $\osched\models \varphi$. 
\end{problem*}
\end{mdframed}
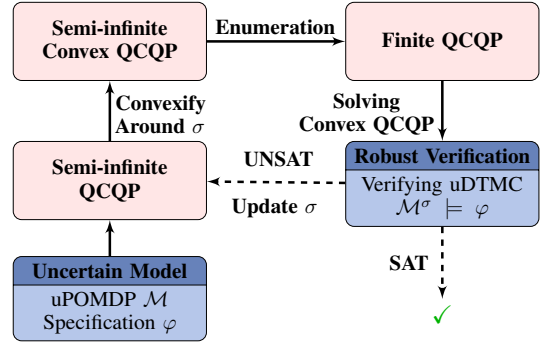
\begin{figure}[t]
		\centering
		\scalebox{0.8}{	\definecolor{bg}{HTML}{ddeedd}
\definecolor{comp}{HTML}{c2d4dd}
\definecolor{impl}{HTML}{b0aac0}
\definecolor{ligb}{HTML}{5E7FC6}
\definecolor{bodybl}{HTML}{85A1DC}
\definecolor{headbl}{HTML}{264C9C}
\definecolor{bgyel}{HTML}{FFDC6B}
\definecolor{bodyyel}{HTML}{FFE58F}
\definecolor{headyel}{HTML}{E9BB25}
\centering
\begin{tikzpicture}[every node/.style={draw, text centered, shape=rectangle, rounded corners, text width=3cm, minimum height=1.3cm, inner sep=3pt}, blank/.style={draw=none,fill=none,ultra thick,inner sep=0.001pt}]
\tikzset{
	splitnode/.style={
		rectangle split,
		rectangle split parts=2,
		rectangle split part fill={headbl!70,bodybl!60}
	}
}
\tikzset{
	normalnode/.style={
		fill=red!10
	}
}
\tikzstyle{split}=[rectangle split,rectangle split parts=2]


\node[normalnode] (infinite-qcqp) {\textbf{Semi-infinite QCQP} 
};

\node[normalnode, above=1.0cm of infinite-qcqp] (infinite-convex-qcqp) {\textbf{Semi-infinite\linebreak Convex QCQP} 
};
\node[normalnode, right=2.3 cm of infinite-convex-qcqp] (finite-qcqp) {\textbf{Finite QCQP}}; 


\node[splitnode, below=1.0 cm of finite-qcqp] (verification) {\textbf{Robust Verification}\nodepart{second} Verifying uDTMC $\pomdp^{\sched} \models \varphi$};

\draw (finite-qcqp) edge[-latex', very thick] (verification);
      \node[blank,above right=-2.55cm and -4.4cm of finite-qcqp] {\textbf{Solving \linebreak Convex QCQP}};

\draw (infinite-qcqp) edge[-latex', very thick] (infinite-convex-qcqp);
\node[blank,above right=-2.55cm and -2.3cm of infinite-convex-qcqp] {\textbf{Convexify \linebreak Around $\sched$}};

\draw (infinite-convex-qcqp) edge[-latex', very thick] (finite-qcqp);

      \node[blank,above right=-1.1cm and -0.45cm of infinite-convex-qcqp] {\textbf{ Enumeration }};




\node[splitnode, below=0.6cm of infinite-qcqp] (pomdp) {\textbf{Uncertain Model} \nodepart{second}uPOMDP $\pomdp$ \\ Specification $\varphi$};

\node[draw=none,below=1.2cm of verification,minimum height=0.1cm,text width =0.5cm] (safe) {\color{green!70!black}\Large\checkmark};


\path [-latex', very thick, dashed] (verification) --node[draw=none,left=0.5cm,text width=0.25cm,minimum height=0.05cm]{\textbf {SAT}} (safe);

\draw (pomdp) edge[-latex', very thick] (infinite-qcqp);

\draw (verification) edge[-latex', very thick, dashed] (safe);
\draw (verification) edge[-latex', very thick, dashed] (infinite-qcqp);

      \node[blank,above right=-1.4cm and -0.4cm of infinite-qcqp] {\textbf{ UNSAT \linebreak\linebreak {Update $\sched$}}};




\end{tikzpicture}}
		\caption{Flowchart of the overall approach.}
		\label{fig:flowchart_problem}
	\end{figure}
\paragraph{Outline.} Fig.~\ref{fig:flowchart_problem} shows the outline of our approach. The input is a uPOMDP $\pomdp$ and one or more specifications $\varphi$. 
We first state a semi-infinite optimization problem which defines the robust synthesis problem within this section.
In Sect.~\ref{sec:ccp}, we show how this nonlinear problem can we convexified around an initial policy $\osched\in\oSched^\pomdp$, followed by Sect.~\ref{sec:uncertainty} which describes how a finite, efficiently solvable problem is obtained. This procedure is augmented by an efficient robust verification method. 

\paragraph{Semi-infinite optimization problem.}
We introduce the following variables for the optimization problem: $\{p_s \mid s \in S\}$ for the probability to reach the targets $T$ from state $s$, $\{c_s \mid s \in S\}$ for the expected cost to reach the goal set $G$ from $s$, and$\{\sched_{s,\alpha} \mid s \in S, \alpha \in Act\}$ to encode the randomized policy.
\begin{flalign}
	&\begin{aligned}
    	\text{minimize} \quad   & f
    \end{aligned} \label{NLP:obj}\\
    &\begin{aligned}
    	\text{subject to} \quad
    	\,\,  p_{s_I} \geq \lambda,  
           \,\,    c_{s_I} \leq \kappa,
    \end{aligned}  \label{NLP:constraint:2}\\
       &\begin{aligned} 
    	\forall s \in T. \quad & p_s = 1,
    	\,\, & \forall s \in G. \quad & c_s = 0, 
    \end{aligned} \label{NLP:constraint:6}\\
    &\begin{aligned} 
    	\forall s \in S. \quad  & \sum\nolimits_{\act \in Act} \sched_{s,\act} = 1,
    \end{aligned} \label{NLP:constraint:3}\\
    &\begin{aligned} 
    	&\forall s \in S \setminus T.\, \forall P \in \mathcal{P}. \\
        & \,\, p_s \leq \sum\nolimits_{\act \in Act} \sched_{s,\act} \cdot \sum\nolimits_{s' \in S} P(s,\alpha,s')\cdot p_{s'},
    \end{aligned} \label{NLP:constraint:5}\\
    &\begin{aligned} &\forall s \in S \setminus G.\, \forall P \in \mathcal{P}. \\
             & \,\, c_s \geq \sum_{\alpha \in Act} \sched_{s,\alpha} \cdot \Big( c(s,\alpha) + \sum_{s' \in S} P(s,\alpha,s') \cdot c_{s'} \Big)
    \end{aligned} \label{NLP:constraint:7} \\
    &\begin{aligned} 
    \forall s, s' \in S.\forall \alpha \in \Act.\ \ObsFun(s) = \ObsFun(s')\rightarrow\sigma_{s,\alpha} = \sched_{s',\alpha}.
    \end{aligned}\label{NLP:constraint:8}
\end{flalign}
Here, $f$ is an objective function, for example $c_{s_I}$ for minimizing the expected cost at the initial state. 
The constraints in \eqref{NLP:constraint:2} encode a reachability and expected cost threshold, respectively. 
\eqref{NLP:constraint:6} defines the fixed reachability and cost values for states that belong to the respective target and goal set, and \eqref{NLP:constraint:3} encodes well-defined policy probabilities.
\eqref{NLP:constraint:5} and \eqref{NLP:constraint:7} define the reachability and cost variables for all other states. 
Note that $p_s$ may at most be assigned the exact probability to reach $T$, ensuring the correct satisfaction of the specification $\reachProplT$.
Finally, \eqref{NLP:constraint:8} defines policy variables from states with the same observation to have the same value, ensuring our policy is based on the observations instead of the states.

We will now take a closer look at the type of this optimization problem.
First, the functions in~\eqref{NLP:constraint:5} and~\eqref{NLP:constraint:7} are \emph{quadratic}. 
Essentially, the policy variables $\sched_{s,\act}$ are multiplied with the probability variables $p_{s}$ in constraint~\eqref{NLP:constraint:5} and with the cost variables $c_{s}$ in \eqref{NLP:constraint:7}.
If we restrict the objective function to be quadratic, we have a QCQP. 
%
As the entries in the transition probability matrices $P(s,\act,s')$ for $s, s' \in S$ and $\act \in \Act$ belong to continuous intervals, there are infinitely many of the constraints~\eqref{NLP:constraint:5} and~\eqref{NLP:constraint:7} over a finite set of variables.
%
%
    Note that we only consider policies where for all states $s$ and actions $\alpha$ it holds that $\sched_{s,\act} > 0$, such that applying the policy to the uPOMDP does not exclude states or transitions.

\section{Convexifying the Semi-Infinite QCQP}
\label{sec:ccp}
We discuss how we convexify the semi-infinite QCQP. 
We use the \emph{penalty convex-concave procedure} (CCP)~\cite{lipp2016variations} which iteratively over-approximates a non-convex optimization problem via linearization.
The resulting convex problem can then be solved efficiently, and the process is iterated until a suitable solution is found.
Specifically, we rewrite the quadratic functions in~\eqref{NLP:constraint:5} and~\eqref{NLP:constraint:7} as a \emph{sum of convex and concave functions} and compute upper bounds for the~\emph{concave} functions.
We check the feasibility regarding the reachability and expected cost specifications using robust value iteration~\cite{benedikt2013ltl}.
Until such a feasible solution is found, the CCP seeks solutions in the vicinity of previous ones.

\begin{remark}
	This section assumes we can effectively solve a semi-infinite convex optimization problem. 
	How to do this in our setting is discussed in Section~\ref{sec:uncertainty}.
\end{remark}

\paragraph{Convexifying the constraints.}
The CCP method starts with any (possibly infeasible) assignment $\hat{p}_s, \hat{c}_s,$ and $\hat{\sched}_{s,\act}$ to the variables $p_s, c_s,$ and $\sched_{s,\act}$.
Consider the bilinear function 
$$h^c(s,\act,s',P)=P(s,\act,s') \cdot \sched_{s,\act}\cdot c_{s'}$$
for any $s, s' \in S, \act \in \Act$ and $P \in \mathcal{P}$ whose right-hand is part of the constraint~\eqref{NLP:constraint:7} in the original QCQP.
For simplicity, we set $P(s,\act,s')=2\cdot d, \sched_{s,\act}=y$, and $c_{s'}=z$ and get $h^c(s,\act,s',P)=2\cdot  d \cdot y \cdot z$.
We rewrite $2\cdot d \cdot y \cdot z$ to 
$2\cdot d \cdot y \cdot z+d(y^2+z^2)-d(y^2+z^2).$
Then, we can write $2\cdot d \cdot y \cdot z+d(y^2+z^2)$ as $h^c_{\textrm{cvx}}(s,\act,s',P)=d(y+z)^2,$ which is a \emph{quadratic convex function}. 
Recalling~\eqref{NLP:constraint:7}, we add the cost function $c(s,\act)$ and get the convex function 
$\hat{h}^c_{\text{cvx}}=h^c_{\text{cvx}}+y \cdot c(s,\act)$
as $y \cdot c(s,\act)$ is affine.

The function $h^c_{\textrm{ccv}}(s,\act,s',P)=-d(y^2+z^2)$, however, is concave, and we have to convexify it.
In particular, we transform $h^c_{\textrm{ccv}}(s,\act,s',P)$ to $\hat{h}^c_{\textrm{ccv}}(s,\act,s',P)=d(\hat{y}^2+\hat{z}^2)+2\cdot d(\hat{y}^2+\hat{z}^2-y\hat{y}-z\hat{z})$, where $\hat{y}$ and $\hat{z}$ are any assignments to the policy and probability variables. 
$\hat{h}^c_{\textrm{ccv}}(s,\act,s',P)$ is affine in $y$ and $z$ and therefore convex.

We convexify~\eqref{NLP:constraint:5} analogously and replace
the quadratic functions with $\hat{h}^p_{\textrm{cvx}}(s,\act,s',P)$ and $\hat{h}^p_{\textrm{ccv}}(s,\act,s',P)$.

After the convexification step, we replace~\eqref{NLP:constraint:5} and~\eqref{NLP:constraint:7} by
\begin{align}
  &\forall s \in S \setminus T.\, \forall P \in \mathcal{P}. \label{Convex:constraint1} \\
&-p_s \geq  \sum_{\alpha \in Act} \sum_{s' \in S}  \big( \hat{h}^p_{\textrm{cvx}}(s,\act,s',P) + \hat{h}^p_{\textrm{ccv}}(s,\act,s',P)\big),\nonumber\\
&\forall s \in S \setminus G.\, \forall P \in \mathcal{P}.  \label{Convex:constraint2}\\
& c_s \geq \sum_{\alpha \in Act} \sum_{s' \in S}  \big( \hat{h}^c_{\textrm{cvx}}(s,\act,s',P) + \hat{h}^c_{\textrm{ccv}}(s,\act,s',P)\big),\nonumber
\end{align}
which are semi-infinite convex constraints in $\sched_{s,\act}$, $p_{s'}$ and $c_{s'}$.
We switch the sign of $p_s$ as it was upper bounded before.

The resulting problem is \emph{convex} (yet semi-infinite). 
As we over-approximate the quadratic functions, any feasible solution to the convex problem is also feasible for the original semi-infinite QCQP.
However, due to the over-approximation, the resulting convex problem might be infeasible though the original one is not.
To find a feasible assignment, we assign a so-called non-negative \emph{penalty variable} $k_{s}$ for each of the probability constraints in~\eqref{Convex:constraint1} for $s \in S\setminus T$, and $l_{s}$ for the cost constraints in~\eqref{Convex:constraint2}. 
To find a solution that induces a minimal infeasibility, or minimal violations to the convexified constraints, we minimize the sum of the penalty variables.
This gives us another semi-infinite convex problem:

\begin{flalign}
&\text{minimize} \,\,\,\,    f +\tau \big(\sum\nolimits_{s \in S\setminus T}k_{s}+\sum\nolimits_{s \in S\setminus G}l_{s}\big) \label{ccp:obj}\\
&\text{subject to}  \quad \,\,  p_{s_I} \geq \lambda,  \,\,    c_{s_I} \leq \kappa,\label{ccp:constraint:1}\\
& \forall s \in T. \quad  p_s = 1,\,\,\,\, \forall s \in G. \quad  c_s = 0,  \label{ccp:constraint:2}\\
&\forall s \in S. \quad   \sum\nolimits_{\alpha \in Act} \sched_{s,\alpha} = 1, \label{ccp:constraint:3}\\
&\forall s \in S \setminus T.\, \forall P \in \mathcal{P}.\label{ccp:constraint:4} \\
& k_{s}-p_s \geq  \sum_{\alpha \in Act} \sum_{s' \in S}  \big( \hat{h}^p_{\textrm{cvx}}(s,\act,s',P) + \hat{h}^p_{\textrm{ccv}}(s,\act,s',P)\big), \nonumber\\
&\forall s \in S \setminus G.\, \forall P \in \mathcal{P}. \label{ccp:constraint:5}\\
&  l_{s}+c_s \geq \sum_{\alpha \in Act} \sum_{s' \in S}  \big( \hat{h}^c_{\textrm{cvx}}(s,\act,s',P) + \hat{h}^c_{\textrm{ccv}}(s,\act,s',P)\big), \nonumber\\
&\forall s, s'\!\in S.\forall \alpha \in \Act.\ \ObsFun(s) = \ObsFun(s')\rightarrow\sigma_{s,\alpha} = \sched_{s',\alpha}\label{ccp:constraint:6}.
\end{flalign}

If a solution assigns all penalty variables to zero, then the solution QCQP is feasible for the original non-convex QCQP, as we over-approximate the concave functions by affine functions.
If any of the penalty variables $k_{s}$ and $l_{s}$ are assigned to a positive value, we update the penalty parameter $\tau$ by $\mu+\tau$ for a $\mu>0$, similar to the approach in~\cite{lipp2016variations}. 
We put an upper limit $\tau_{\text{max}}$ on $\tau$ to avoid numerical problems during the procedure. 
After getting a new assignment, we convexify the non-convex QCQP by linearizing the concave functions around the new assignment, and solve the resulting convex QCQP.
We repeat the procedure until we find a feasible solution.
If the CCP converges to an infeasible solution, we restart the procedure with another value of the policy $\hat{\sched}$.
Note that convergence is guaranteed for a fixed $\tau$, i.e, after $\tau=\tau_{\text{max}}$, but it may converge to an infeasible point of the original problem~\cite{lipp2016variations}.

\section{Derivation of the Convex QCQP}\label{sec:uncertainty}%
%
In this section, we describe how to transform the semi-infinite convex QCQP to a convex QCQP that is amenable to efficient solving techniques. 
That transformation largely depends on the type of uncertainty set that enters the problem.
So far, we just stated that we have to account for all concrete probability functions $P$ within the uncertainty set $\mathcal{P}$, see the constraints~\eqref{NLP:constraint:5} and~\eqref{NLP:constraint:7}.
We now describe this abstract notion as a specific uncertainty set.
In particular, we enumerate exactly all possible (valid) probability distributions from the uncertainty set that enter the problem in the form of probability intervals. 


For a uPOMDP $\PomdpInit$ and its underlying uMDP $\umdpInit$, each state-action pair has a fixed number of associated probability intervals.
For state $s\in S$ and action $\act\in\Act$, we assume $n$ intervals $[a_i, b_i]\in\intervals$, $1\leq i\leq n$.

For each transition probability function $P\in\mathcal{P}$ at $(s,\act)$, we ensure that $P$ is valid via 
    $\forall P \in \mathcal{P}. \quad \sum\nolimits_{s' \in S} P(s,\alpha,s') = 1$. 



We define the set of all possible probability distributions formed by the intervals $[a_i, b_i]$, expressed by the following set of linear constraints that form a \emph{convex polytope}:
\begin{align}
    \forall i, 1\leq i\leq n. \quad  a_i \leq x_i \leq b_i, \,\,\,\, \sum\nolimits_{i = 1}^n x_i = 1.
\end{align}
We rewrite these constraints into their canonical form of $A\vec{x} \leq \vec{c}$. 
Note that
$a_i \leq x_i \leq b_i$ can be rewritten as
\begin{align}
    \forall i. \,\, -x_i \leq -a_i,\, \,\,\forall i.\,\,  x_i \leq b_i.
\end{align}
The equality constraint can be rewritten as the conjunction of $\leq$ and $\geq$. 
Finally, multiplying by $-1$ will flip the $\geq$ sign:
\begin{align}
    \quad \sum\nolimits_{i = 1}^n x_i \leq 1, \, \,\,\sum\nolimits_{i = 1}^n -x_i \leq -1.
\end{align}
From these constraints we construct matrix $A$ and vector $\vec{c}$:
\begin{align}
  &  A^{\top} = \begin{bmatrix}
        -I_n & I_n & H_n^{\top} & -H_n^{\top}
    \end{bmatrix},\\
   & \vec{c}^{\top} = \begin{bmatrix}
        -a_1 &
        \cdots &
        -a_n &
        b_1 &
        \cdots &
        b_n &
        1 &
        -1
    \end{bmatrix}
\end{align}
$I_n$ is the $n \times n$ identity matrix, and $H_n$ is the $1 \times n$ single row matrix consisting of only ones. 
This matrix and vector are the canonical form to describe a convex polytope.
We enumerate all the vertices of this convex polytope by using the double description method~\cite{Fukuda1996doubledescr}.

\paragraph{Exact representation of distributions.}
By construction, the polytope describes \emph{exactly the set of valid probability distributions} from the intervals.
Moreover, because the polytope is convex, it suffices to enumerate over the vertices \cite{Lofberg2012automaticrobust} to capture all of these distributions. 
As we have these vertices now, we can simply replace the robust constraints \eqref{ccp:constraint:4} and \eqref{ccp:constraint:5} by a (finite) number of constraints in which the uncertainty is replaced by all possible combinations of values from the vertices.
Effectively we enumerate all possible probabilities against which the policy needs to be robust. 
The resulting convex QCQP can directly be solved, for instance, by the QCQP solver \tool{Gurobi}~\cite{gurobi}.

\paragraph{Complexity of the convex QCQP.} Note that solving a robust convex QCQP with polytopic uncertainty is still NP-Hard~\cite{bertsimas2011theory} as the number of vertices of a convex polytope can be exponential in the number of dimensions.
However, in our specific case where we apply this method to uPOMDPs, the dimension of each polytope is determined by the number of successor states for each state-action pair. 
As mentioned before, we expect the number of successors to be low, and thus the dimension of each polytope and the number of vertices, to be manageable.

\paragraph{Integrating robust verification to the penalty CCP.}
In each iteration of the CCP as described in Sect.\ref{sec:ccp}, the QCQP solver assigns concrete values to the variables which induces a concrete policy $\osched\in\oSched^\pomdp$. 
We apply this instantiation to the uPOMDP $\pomdp$ and an \emph{uncertain Markov chain} $\pomdp^\osched$ results.
For this model without partial observability and nondeterminism, we employ \emph{robust value iteration}~\cite{wiesemann2013robust} to check whether the specifications are already satisfied.
Our numerical examples show that this additional verification step is a good heuristic for an earlier termination of the CCP when the penalty variables have not evaluated to zero yet.
We also use the result to ensure that the probability and the cost variables are consistent with the policy variables for the next iteration of the CCP.

\section{Numerical Examples}\label{sec:experiments}
\begin{table*}[t]
	\setlength{\tabcolsep}{2pt}
	\caption{Numerical examples.}
	\label{tab:ToolComp}
	\centering
	\scalebox{0.92}{%
		\begin{tabular}{@{}cc|rrrr|rrr|rrr@{}}
			\toprule
			& & \multicolumn{4}{c}{Nominal}   &\multicolumn{3}{|c}{Small Interval} & \multicolumn{3}{|c}{Big Interval}  \\
			Problem      & Type   & States & Constraints & Iteration & Time (s)& Constraints & Iteration & Time (s) & Constraints & Iteration & Time (s) \\
			\midrule
			Maze &$\Ex_{\leq 80}$  & 30 &    27 & 40 &   1.46 & 68 & 100 & 4.68 & 68 & 275 & 13.90  \\
			Maze &$\Ex_{\leq 50}$  & 30 &    27 & 42 &   1.84 & 68 & 101 & 4.90 & 68 & 1222 & 52.07  \\
			Maze &$\Ex_{\leq 25}$  & 30 &    27 & 49 &   1.35 & 68 & 104 & 4.67 & 68 & 2122 & 105.54  \\
			Grid &$\mathbb{P}_{\geq 0.84}$  & 18  & 14 &  8 & 0.11 &   56 & 8 & 0.20 & 56 & 23 & 2.15  \\			
			Grid &$\mathbb{P}_{\geq 0.92}$  & 18 &   14 & 8 & 0.11 &  56 & 9 & 0.22 & 56 & 56 & 6.26  \\
			Aircraft &$\mathbb{P}_{\geq 0.80}$  & 175861 & 214448 & 2 & 5.89 & 399094 & 5 & 39.61 & 399004 & 13 & 106.72  \\
						Aircraft &$\mathbb{P}_{\geq 0.90}$  & 175861 & 214448 & 5 & 31.98 & 399004 & 26 & 215.94 & 399004 & 59 & 540.33  \\
				Aircraft &$\mathbb{P}_{\geq 0.93}$  & 175861 & 214448 & 22 & 136.64 & 399004 & 67 & 637.23 & 399004 & 143 & 1314.29  \\
				Aircraft &$\mathbb{P}_{\geq 0.95}$  & 175861 & 214448 & 40 & 274.43 & 399004 & 172 & 1475.70 & 399004 & TO & TO  \\
			Aircraft &$\mathbb{P}_{\geq 0.97}$  & 175861 & 214448 & 50 & 339.78 & 399004 & TO & TO & 399004 & TO & TO  \\	
				Network &$\mathbb{E}_{\leq 90}$  & 38719 & 107068  &  8 &   100.57 & 187591 & 8 & 114.99 & 187591  & 8 & 234.57  \\
			Network &$\mathbb{E}_{\leq 50}$  & 38719 & 107068  &  10 &   118.51 & 187591 & 10 & 148.98 & 185791 & 10 & 291.18  \\
						Network &$\mathbb{E}_{\leq 5}$  & 38719 & 107068   &  12 &   135.25 & 187591 & 12 & 171.78 & 187591  & 12 & 336.33  \\
			\bottomrule
	\end{tabular}}
\end{table*}
We evaluate our robust synthesis procedure on benchmark examples that are subject to either reachability or expected cost specifications.
As part of a \tool{Python} toolchain, we use the probabilistic model checker \tool{Storm}~\cite{DBLP:conf/cav/DehnertJK017} to extract an explicit state space representation of uPOMDPs.
The experiments were performed on a computer with an Intel Core i9-9900u 2.50 GHz processor and 64 GB of RAM with Gurobi 9.0~\cite{gurobi} as the QCQP solver and our own implementation of a robust value iteration. 
We use a 1 hour time-out (TO).
For all our examples we use a standard POMDP model so-called \emph{nominal} probabilities, as well as two different sizes of probability intervals, namely a small one and a big one. 
The reason for these three options is that we want to showcase the effect of growing uncertainty on the runtime governed by the number of CCP iterations.

\begin{figure}
\centering
	\subfigure[Grid-world.]
	{\scalebox{0.43}{
		\input{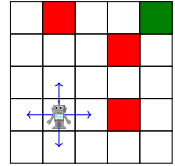}
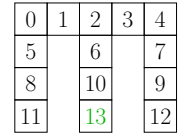\label{fig:gridworld}
	}}\hspace{2cm}
	\subfigure[Maze.]
	{\scalebox{0.43}{
		\begin{tikzpicture}
    [
        box/.style={rectangle,draw=black,thick, minimum size=1cm},
    ]
\colorlet{darkgreen}{green!70!black}

        \node[box] at (0,0){};
        \node[box] at (1,0){};
        \node[box] at (2,0){};
        \node[box] at (3,0){};
        \node[box] at (4,0){};
        \node[box] at (4,-1){};
        \node[box] at (4,-2){};
        \node[box] at (4,-3){};
        \node[box] at (2,-1){};
        \node[box] at (2,-2){};
        \node[box] at (2,-3){};

        \node[box] at (0,-1){};
        \node[box] at (0,-2){};
        \node[box] at (0,-3){};

  \node [draw=none] at (0,0) {\Huge $0$};
  \node [draw=none] at (1,0) {\Huge $1$};
  \node [draw=none] at (2,0) {\Huge $2$};
  \node [draw=none] at (3,0) {\Huge $3$};
  \node [draw=none] at (4,0) {\Huge $4$};
    \node [draw=none] at (0,-1) {\Huge $5$};
    \node [draw=none] at (2,-1) {\Huge $6$};
    \node [draw=none] at (4,-1) {\Huge $7$};
    \node [draw=none] at (0,-2) {\Huge $8$};
    \node [draw=none] at (4,-2) {\Huge $9$};
    \node [draw=none] at (2,-2) {\Huge ${10}$};
    \node [draw=none] at (0,-3) {\Huge $11$};
    \node [draw=none] at (4,-3) {\Huge $12$};
    \node [draw=none] at (2,-3) {\Huge ${\color{darkgreen}13}$};

\end{tikzpicture}\label{fig:maze}
	}}
\label{fig:examples}
\caption{Two standard POMDP examples.}
\end{figure}

\paragraph{Standard POMDP examples.}
We consider the following POMDP case studies with added uncertainties in the form of intervals.
\emph{Grid-world robot} is based on the POMDP example in~\cite{littman1995learning}.
A robot is placed randomly into a grid and the goal is to safely reach the north-east corner, see Fig.~\ref{fig:gridworld}. 
The robot may only reach intended states with a certain probability. 
We consider three variants for that probability: $0.98$ (nominal), $[0.95,0.98]$ (Small Interval), and $[0.50,0.98]$ (Big Interval), yielding two distinct uPOMDPs and one POMDP.
The reachability specification $\mathbb{P}_{\geq 
\lambda}$ ensures to reach the target without visiting the ``traps'.

%
%
%

The second example is a maze setting, introduced in~\cite{mccallum1993overcoming}, where a robot is to reach a target location in minimal time see Fig.~\ref{fig:maze}.
Again, we consider a ``slippery'' maze, similar to the previous example.
We use the following probabilities for slipping: $0.97$ (nominal), $[0.94,0.97]$ (small interval), and $[0.50,0.97]$ (big interval).
We define an expected cost specification $\Ex_{\leq \kappa}$ to reach the goal.

Our third example is from~\cite{yang2011real} and concerns scheduling wireless traffic, where at each time period a scheduler generates a new packet for each user. 
We assume the scheduler does not know the exact probabilities of the channels between the scheduler and the users, therefore, the transition probabilities between the channel states belong to intervals.
The scheduler does not now the current states of the users, and has to schedule wireless traffic based on the partial information of the states of the users with $9743$ possible observations.
The nominal probability is $0.9$, and we have $[0.875,0.9]$ (Small Interval) and $[0.8,0.9]$ (Big Interval).
The specification $\mathbb{E}_{\leq \kappa}$ is to minimize the expected number of dropped packets for the scheduler.

\paragraph{Aircraft collision avoidance.}
In this more sophisticated example, we consider a robust aircraft collision avoidance problem~\cite{kochenderfer2015decision}.
The objective is to maximize the probability of avoiding a collision with an intruder aircraft while taking into account sensor errors and uncertainty in the future paths of the intruder.
The problem is a POMDP with state variables (1) $h$, altitude of the intruder relative to the own aircraft, (2) $\dot{h}$, vertical rate of the intruder relative to the own aircraft, (3) $\tau$, time to potential collision, (4) $s_{\textrm{res}}$, whether the pilot is responsive to requested commands.

We discretize the $h$ variable into 33 points over the range $\pm 4000$ feet, the $\dot{h}$ variable into 25 points between $\pm 10,000$ feet/minute and $\tau$ to 40 points from 0 to 40 seconds.
The $1905$ possible observations give partial information of $h$ and $\dot{h}$.
In the POMDP model, the probability of getting a correct observation  is $0.95$.
Again, we assess the effect of the interval size by means of two intervals, namely $[0.90,0.95]$ (Small Interval) and $[0.75,0.95]$ (Big Interval).
The specification $\mathbb{P}_{\geq \lambda}$ is to maximize the probability of not having a collusion with the intruder within $40$ seconds. 
Similarly, we use different values of $\lambda$ to show the effect of different probability thresholds.


\paragraph{Discussion of the results.}
In Table~1, we list the experimental results for different specification thresholds for each example. 
\textquotedblleft States\textquotedblright ~denotes the number of states in the model, \textquotedblleft Constraints\textquotedblright ~denotes the number of constraints in the convex QCQP, \textquotedblleft Iterations\textquotedblright ~denotes the number of CCP iterations, and \textquotedblleft Time (s)\textquotedblright~denotes the time spent in Gurobi in seconds.
We pick the specification thresholds such that one is near to the point where our procedure converges to an infeasible solution. 

We remark that the number of constraints in each example increases by adding intervals (instead of concrete probabilities) to the model, due to the explicit enumeration of polytope vertices, see Section~\ref{sec:uncertainty}.
However, the number of constraints does not depend on the size of the intervals.
We also note that the solution time for each iteration for the problem with uncertainty (uPOMDP) is larger than for the original model (POMDP) due to these additional constraints.

For the examples with small state spaces, namely Maze and Grid, we picked the thresholds $25$ and $0.92$, respectively, to be very near the threshold where our method converges to an infeasible solution for the case with a larger uncertainty. 
We observe that the number of iterations may grow rapidly with a decreasing threshold.
In particular, already for a threshold $25$ for the Maze example, the number of iterations for the case with a larger uncertainty is much bigger compared to the nominal and the small uncertainty case.


For the Aircraft example, we observe that the number of iterations required to satisfy a reachability specification increases significantly with an increasing degree of uncertainty. 
For threshold $0.95$, we cannot find a policy that induces a reachability probability that is larger than the threshold with large uncertainty. 
Similarly, for $0.97$, both of the uPOMDPs converged to a policy that does not satisfy the specification.

On the other hand, for the Network example, the number of iterations for all three cases is the same with three different models, and the only difference between the cases is the computation time per iteration.
Particularly, the optimization problems with larger uncertainty were more numerically challenging for Gurobi, and computing the optimal solution for the optimization problems took more time per iteration.

%
%
%



\section{Conclusion and Future Work}\label{sec:conclusion}
We presented a new approach to computing robust policies for uncertain POMDPs. 
The experiments showed that we are able to apply our method based on convex optimization on well-known benchmarks with varying levels of uncertainty. 

Future work will incorporate finite memory into the robust policies, similar to~\cite{junges2018finite}.
Basically, a memory structure can be directly integrated into the state space of the uPOMDP rendering our method applicable without further adaption.
We will also combine our method with model-based reinforcement learning where the uncertainty accounts for insufficient confidence in transition probabilities.

%
\begin{small}
\bibliographystyle{aaai}
\bibliography{literature}
\end{small}


\end{document}